\documentclass{amsart}
\usepackage{amsmath}
\usepackage{eurosym}
\usepackage{amsfonts}
\usepackage{amssymb}
\usepackage{amsmath}
\usepackage{amsfonts}
\usepackage{graphicx}
\usepackage{subfig}
\usepackage{eurosym}
\usepackage{amssymb}
\usepackage{amsmath}
\usepackage{amsfonts}
\usepackage{graphicx}
\usepackage{subfig}

\newtheorem{theorem}{Theorem}
\theoremstyle{plain}

\newtheorem{definition}{Definition}

\numberwithin{equation}{section}

\begin{document}
\title[]{The Clairaut's theorem on rotational surfaces in pseudo Euclidean
4-space with index 2}
\author{Fatma ALMAZ}
\address{Department of Mathematics, Firat university, 23119 ELAZI\u{G}/T\"{U}%
RK\.{I}YE}
\email{fb\_fat\_almaz@hotmail.com}
\author{M\.{I}hr\.{I}ban ALYAMA\c{C} K\"{U}LAHCI}
\address{Department of Mathematics, Firat university, 23119 ELAZI\u{G}/T\"{U}%
RK\.{I}YE}
\email{mihribankulahci@gmail.com }
\thanks{This paper is in final form and no version of it will be submitted
for publication elsewhere.}

\begin{abstract}

In this paper, Clairaut's theorem is expressed on the
surfaces of rotation in semi Euclidean 4-space. Moreover, the general
equations of time-like geodesic curves are characterized according to the
results of Clairaut's theorem on the hyperbolic surfaces of rotation and the
elliptic surface of rotation, respectively.
\end{abstract}

\maketitle

\section{Introduction}

The geodesics for rotational surfaces have been studied for a long time and
many examples of rotational surfaces have been discovered. To understand the
rest of geodesics; we need Clairaut's Theorem, which is
very helpful to understand the geodesics on surfaces of rotation. This gives
a well-known characterization of geodesics on surfaces of rotation.

Many studies of surfaces of rotation have received much attention from our
researchers. Among them, one can cite our work \cite{1}, we described the
rotational surfaces using curves and matrices which are the subgroups of
rotating a selected axis in Galilean 4-space. We examined the tube surfaces
generated by the curve in Galilean 3-space and gave certain results of
describing the geodesics on the surfaces \cite{2,4}. We gave the surfaces of
rotational generated by a magnetic curve. Also, we gave the conditions being
geodesic on these rotational surfaces in null cone 3-space, with the help of
Clairaut's theorem \cite{3}. In our study \cite{5} we expressed the
hyperbolic and the elliptic rotational surfaces using a curve and matrices
in 4-dimensional semi-Euclidean space. Goemans constructed a new type of
surfaces in Euclidean and Lorentz\^{a}\euro ``Minkowski 4-space and proved
the classification theorems of flat double rotational surfaces \cite{8}.
Hoffmann and Zhou \cite{9}, discussed some issues of displaying 2D surfaces
in 4D space, including the behaviour of surface normals under projection,
the silhouette points due to the projection, and methods for object
orientation and projection center specification.

\section{Preliminaries}

Let $E_{2}^{4}$ denote the $4-$dimensional pseudo-Euclidean space with
signature $(2,4)$, that is, the real vector space $%
\mathbb{R}
^{4}$ endowed with the metric $\left\langle ,\right\rangle _{E_{2}^{4}}$
which is defined by 
\begin{equation}
\left\langle ,\right\rangle
_{E_{2}^{4}}=-dx_{1}^{2}-dx_{2}^{2}+dx_{3}^{2}+dx_{4}^{2}  \tag{2.1}
\end{equation}%
or 
\begin{equation}
g=%
\begin{bmatrix}
-1 & 0 & 0 & 0 \\ 
0 & -1 & 0 & 0 \\ 
0 & 0 & 1 & 0 \\ 
0 & 0 & 0 & 1%
\end{bmatrix}%
,  \tag{2.2}
\end{equation}%
where $(x_{1},x_{2},x_{3},x_{4})$ is a standard rectangular coordinate
system in $E_{2}^{4}$.

Recall that an arbitrary vector $v\in E_{2}^{4}\backslash \{0\}$ can have
one of three characters: it can be space-like if $g(v,v)>0$ or $v=0,$
time-like if $g(v,v)<0$ and null if $g(v,v)=0$ and $v\neq 0.$

The norm of a vector $v$ is given by $\parallel v\parallel =\sqrt{g(v,v)}$
and two vectors $v$ and $w$ are said to be orthogonal if $g(v,w)=0$. A
space-like or time-like curve $x(s)$ has unit speed, if $g(x^{\prime},x^{%
\prime })=\pm 1.$

Let $%
(x_{1},x_{2},x_{3},x_{4}),(y_{1},y_{2},y_{3},y_{4}),(z_{1},z_{2},z_{3},z_{4}) 
$ be any three vectors in $E_{2}^{4}$. The pseudo Euclidean cross product is
given as 
\begin{equation*}
x\wedge y\wedge z=%
\begin{pmatrix}
-i_{1} & -i_{2} & i_{3} & i_{4} \\ 
x_{1} & x_{2} & x_{3} & x_{4} \\ 
y_{1} & y_{2} & y_{3} & y_{4} \\ 
z_{1} & z_{2} & z_{3} & z_{4}%
\end{pmatrix}%
,
\end{equation*}%
where $i_{1}=\left( 1,0,0,0\right) ,i_{2}=\left( 0,1,0,0\right)
,i_{3}=\left( 0,0,1,0\right) ,i_{4}=\left( 0,0,0,1\right) $, \cite{9}.

The pseudo-Riemannian sphere $S_{2}^{3}\left( m,r\right) $ centred at $m\in
E_{2}^{4}$ with radius $r>0$ of $E_{2}^{4}$ is defined by 
\begin{equation*}
S_{2}^{3}\left( m,r\right) =\left\{ x\in E_{2}^{4}:\left\langle
x-m,x-m\right\rangle =r^{2}\right\} .
\end{equation*}

The pseudo-hyperbolic space $H_{1}^{3}\left( m,r\right) $ centred at $m\in
E_{2}^{4}$ with radius $r>0$ of $E_{2}^{4}$ is defined by%
\begin{equation*}
H_{1}^{3}\left( m,r\right) =\left\{ x\in E_{2}^{4}:\left\langle
x-m,x-m\right\rangle =-r^{2}\right\} .
\end{equation*}

The pseudo-Riemannian sphere $S_{2}^{3}\left( m,r\right) $ is diffeomorphic
to $%
\mathbb{R}
^{2}\times S$ and the pseudo-hyperbolic space $H_{1}^{3}\left( m,r\right) $
is diffeomorphic to $S^{1}\times 
\mathbb{R}
^{2}$. The hyperbolic space $H^{3}\left( m,r\right) $ is given by%
\begin{equation*}
H^{3}\left( m,r\right) =\left\{ x\in E_{2}^{4}:\left\langle
x-m,x-m\right\rangle =-r^{2},x_{1}>0\right\},
\end{equation*}
\cite{7,12,14}.


\begin{definition}
\cite{10}, A one-parameter group of diffeomorphisms of a manifold $M$ is a
regular map $\psi:M\times%
\mathbb{R}
\rightarrow M$, such that $\psi_{t}(x)=\psi(x,t),$ where

\begin{enumerate}
\item $\psi_{t}:M\rightarrow M$ is a diffeomorphism

\item $\psi_{0}=id$

\item $\psi_{s+t}=\psi_{s}o\psi_{t}.$
\end{enumerate}
\end{definition}

This group is attached with a vector field $W$ given by $\frac{d}{dt}\psi
_{t}(x)=W(x),$ and the group of diffeomorphism is said to be the flow of $W$.

\begin{definition}
If a one-parameter group of isometries is generated by a vector field $W$,
then this vector field is called as a Killing vector field, \cite{10}.
\end{definition}

\begin{definition}
Let $W$ be a vector field on a smooth manifold $M$ and $\psi_{t}$ be the
local flow generated by $W$. For each $t\in%
\mathbb{R}
,$ the map $\psi_{t}$ is a diffeomorphism of $M$ and given a function $f$ on 
$M $, one considers the Pull-back $\psi_{t}f$ and one defines the Lie
derivative of the function $f$ as to $W$ by%
\begin{equation}
L_{_{W}}f=\underset{t\longrightarrow0}{\lim}\underset{}{\left( \frac{\psi
_{t}f-f}{t}\right) =\frac{d\psi_{t}f}{dt}_{t=0}}.  \tag{2.4}
\end{equation}

Let $g_{\xi\varrho}$ be any pseudo-Riemannian metric, then the derivative is
given as 
\begin{equation*}
L_{_{W}}g_{\xi\varrho}=g_{\xi\varrho,z}W^{z}+g_{\xi
z}W_{,\varrho}^{z}+g_{z\varrho}W_{,\xi}^{z}.
\end{equation*}

In Cartesian coordinates in Euclidean spaces where $g_{\xi\varrho,z}=0,$ and
the Lie derivative is given by%
\begin{equation*}
L_{_{W}}g_{\xi\varrho}=g_{\xi z}W_{,\varrho}^{z}+g_{z\varrho}W_{,\xi}^{z}.
\end{equation*}

In \cite{6,10,11,15}, the vector $W$ generates a Killing field if and if
only 
\begin{equation*}
L_{_{W}}g=0.
\end{equation*}
\end{definition}

\begin{theorem}
Let $\xi $ be a geodesic on a surface of revolution $\Upsilon$ and let $\rho 
$ be the distance function of a point of $\Upsilon$ from the axis of
rotation, and let $\theta$ be the angle between $\xi $ and the meridians of $%
\Upsilon$. Then, $\rho sin\theta $ is constant along $\xi $. Conversely, if $%
\rho sin\theta $ is constant along some curve $\xi $ in the surface, and if
no part of $\xi $ is part of some parallel of $\Upsilon$, then $\xi $ is a
geodesic, \cite{13}.
\end{theorem}

\begin{theorem}
Let the pseudo Euclidean group be a subgroup of the diffeomorphisms group in 
$E_{2}^{4}$ and let $W$ be vector field which generate the isometries. Then,
the killing vector field associated with the metric $g$ is given as%
\begin{equation*}
W(\xi ,\varrho ,\vartheta ,\eta )=a\left( \eta \partial \xi +\xi \partial
\eta \right) +b\left( \vartheta \partial \varrho +\varrho \partial \vartheta
\right) +c\left( \vartheta \partial \xi +\xi \partial \vartheta \right)
\end{equation*}%
\begin{equation*}
+d(\eta \partial \varrho +\varrho \partial \eta )+e(\vartheta \partial \eta
-\eta \partial \vartheta )+f\left( \xi \partial \varrho -\varrho \partial
\xi \right) ,
\end{equation*}%
where $a,b,c,d,e,f\in 
\mathbb{R}
_{0}^{+},$ \cite{5}.
\end{theorem}

\begin{theorem}
Let $W(\xi ,\varrho ,\vartheta ,\eta )$ be the killing vector field and let $%
\gamma=(f_{1},f_{2},f_{3},f_{4})$ be a curve in $E_{2}^{4}$, then the
surfaces of rotation are given as follows

\begin{enumerate}
\item For the rotations $\Omega _{1}=\vartheta \partial \xi +\xi \partial
\vartheta $ and $\Omega _{4}=\eta \partial \varrho +\varrho \partial \eta ,$
the hyperbolic surface of rotation is given as%
\begin{equation*}
S_{14}(x,\alpha ,s)=\left( 
\begin{array}{c}
f_{1}\cosh x+f_{3}\sinh x,f_{2}\cosh \alpha +f_{4}\sinh \alpha , \\ 
f_{1}\sinh x+f_{3}\cosh x,f_{2}\sinh \alpha +f_{4}\cosh \alpha 
\end{array}%
\right) 
\end{equation*}%
and for the planar curve $\gamma (s)=(f_{1}(s),0,0,f_{4}(s))$ the Gaussian
curvature $K$ and the mean curvature vector $H$ of the rotational surface $%
S_{14}(x(t),\alpha (t),s)=\left( f_{1}\cosh x,f_{4}\sinh \alpha ,f_{1}\sinh
x,f_{4}\cosh \alpha \right) $ are given as%
\begin{equation*}
K=\frac{\left( f_{1}^{\prime }f_{4}^{{}}-f_{1}f_{4}^{\prime }\right)
^{2}\left( \overset{.}{x}\overset{.}{\alpha }\right) ^{2}}{f_{4}^{2}\overset{%
.}{\alpha }^{2}-f_{1}^{2}\overset{.}{x}^{2}}+\frac{\left( f_{1}^{\prime
}f_{4}\overset{.}{\alpha }^{2}-f_{4}^{\prime }f_{1}\overset{.}{x}^{2}\right)
\left( f_{1}^{\prime }f_{4}^{\prime \prime }-f_{1}^{\prime \prime
}f_{4}^{\prime }\right) }{-f_{1}^{\prime 2}+f_{4}^{\prime 2}},
\end{equation*}%
\begin{equation*}
H=\{\frac{f_{1}f_{4}\left( \overset{..}{x}\overset{.}{\alpha }+\overset{.}{x}%
\overset{..}{\alpha }\right) }{2\sqrt{f_{4}^{2}\overset{.}{\alpha }%
^{2}-f_{1}^{2}\overset{.}{x}^{2}}}+\frac{f_{4}^{\prime }f_{1}\overset{.}{x}%
^{2}-f_{1}^{\prime }f_{4}\overset{.}{\alpha }^{2}}{2\sqrt{-f_{1}^{\prime
2}+f_{4}^{\prime 2}}}\}e_{3}+\frac{\left( f_{1}^{\prime }f_{4}^{\prime
\prime }-f_{1}^{\prime \prime }f_{4}^{\prime }\right) }{2\sqrt{%
-f_{1}^{\prime 2}+f_{4}^{\prime 2}}}e_{4}
\end{equation*}%
where 

$e_{3}=\frac{\left( f_{4}\overset{.}{\alpha }\sinh x,f_{1}\overset{.}{x%
}\cosh \alpha ,f_{4}\overset{.}{\alpha }\cosh x,f_{1}\overset{.}{x}\sinh
\alpha \right) }{\sqrt{f_{4}^{2}\overset{.}{\alpha }^{2}-f_{1}^{2}\overset{.}%
{x}^{2}}},e_{4}=\frac{\left( f_{4}^{\prime }\cosh x,f_{1}^{\prime }\sinh
\alpha ,f_{4}^{\prime }\sinh x,f_{1}^{\prime }\cosh \alpha \right) }{\sqrt{%
-f_{1}^{\prime 2}+f_{4}^{\prime 2}}}.$

\item For the rotations $\Omega _{2}=\eta \partial \xi +\xi \partial \eta $
and $\Omega _{3}=\vartheta \partial \varrho +\varrho \partial \vartheta ,$
the hyperbolic surface of rotation is given as 
\begin{equation*}
S_{23}(y,z,s)=\left( 
\begin{array}{c}
f_{1}\cosh y+f_{4}\sinh y,f_{2}\cosh z+f_{3}\sinh z, \\ 
f_{2}\sinh z+f_{3}\cosh z,f_{1}\sinh y+f_{4}\cosh y%
\end{array}%
\right) .
\end{equation*}%
and for the planar curve $\gamma (s)=(f_{1}(s),f_{2}(s),0,0)$ the Gaussian
curvature $K$ and the mean curvature vector $H$ of the rotational surface $%
S_{23}(y(t),z(t),s)=\left( f_{1}\cosh y,f_{2}\cosh z,f_{2}\sinh z,f_{1}\sinh
y\right) $ are given as%
\begin{equation*}
K=-\left( 
\begin{array}{c}
\frac{\left( f_{1}f_{2}^{\prime }+f_{1}^{\prime }f_{2}^{{}}\right)
^{2}\left( \overset{.}{y}\overset{.}{z}\right) ^{2}}{f_{2}^{2}\overset{.}{z}^{2}%
+f_{1}^{2}\overset{.}{y}^{2}}+ \\ 
\frac{\left( f_{1}^{{}}f_{2}^{\prime }\overset{.}{y}^{2}+f_{1}^{\prime }f_{2}%
\overset{.}{z}^{2}\right) \left( f_{1}^{\prime \prime }f_{2}^{\prime
}+f_{1}^{\prime }f_{2}^{\prime \prime }\right) }{f_{1}^{\prime
2}+f_{2}^{\prime 2}}%
\end{array}%
\right) ;H=\left( 
\begin{array}{c}
\frac{f_{1}^{{}}f_{2}(\overset{.}{y}\overset{..}{z}+\overset{..}{y}\overset{.%
}{z})}{2\sqrt{f_{2}^{2}\overset{.}{z}^{2}+f_{1}^{2}\overset{.}{y}^{2}}}e_{3} \\ 
+\frac{f_{1}^{{}}f_{2}^{\prime }\overset{.}{y}^{2}+f_{1}^{\prime }f_{2}%
\overset{.}{z}^{2}-f_{1}^{\prime \prime }f_{2}^{\prime }-f_{1}^{\prime
}f_{2}^{\prime \prime }}{2\sqrt{f_{1}^{\prime 2}+f_{2}^{\prime 2}}}e_{4}%
\end{array}%
\right) ,
\end{equation*}%
where 

$e_{3}=\frac{\left( f_{2}\overset{.}{z}\sinh y,f_{1}\overset{.}{y}%
\sinh z,f_{1}\overset{.}{y}\cosh z,f_{2}\overset{.}{z}\cosh y\right) }{\sqrt{%
f_{2}^{2}\overset{.}{z}^{2}+f_{1}^{2}\overset{.}{y}^{2}}}$; $e_{4}=\frac{\left(
f_{2}^{\prime }\cosh y,f_{1}^{\prime }\cosh z,f_{1}^{\prime }\sinh
z,f_{2}^{\prime }\sinh y\right) }{\sqrt{f_{1}^{\prime 2}+f_{2}^{\prime 2}}}$

\item For the rotations $\Omega _{5}=\xi \partial \varrho -\varrho \partial
\xi $ and $\Omega _{6}=\vartheta \partial \eta -\eta \partial \vartheta ,$
the elliptic surface of rotation is given as 
\begin{equation*}
S_{56}(\beta ,\theta ,s)=\left( 
\begin{array}{c}
f_{1}\cos \beta +f_{2}\sin \beta ,-f_{1}\sin \beta +f_{2}\cos \beta , \\ 
f_{3}\cos \theta +f_{4}\sin \theta ,-f_{3}\sin \theta +f_{4}\cos \theta 
\end{array}%
\right) ,
\end{equation*}%
and for the planar curve $\gamma (s)=(0,f_{2}(s),0,f_{4}(s))$ the Gaussian
curvature $K$ and the mean curvature vector $H$ of the rotational surface $%
S_{56}(\beta \left( t\right) ,\theta \left( t\right) ,s)=\left( f_{2}\sin
\beta ,f_{2}\cos \beta ,f_{4}\sin \theta ,f_{4}\cos \theta \right) $ are
given as 
\begin{equation*}
K=-\left( 
\begin{array}{c}
\frac{\left( f_{2}^{\prime }f_{4}-f_{2}f_{4}^{\prime }\right) ^{2}\left( 
\overset{.}{\beta }\overset{.}{\theta }\right) ^{2}}{-f_{2}^{2}\overset{.}{%
\beta }^{2}+f_{4}^{2}\overset{.}{\theta }^{2}}+
\frac{\left( -f_{2}^{\prime \prime }f_{4}^{\prime }+f_{2}^{\prime
}f_{4}^{\prime \prime }\right) (f_{4}^{\prime }f_{2}\overset{.}{\beta }%
^{2}-f_{2}^{\prime }f_{4}\overset{.}{\theta }^{2})^{2}}{-f_{2}^{\prime
2}+f_{4}^{\prime 2}}%
\end{array}%
\right) ;
\end{equation*}
\begin{equation*}
H=\frac{f_{4}f_{2}\left( \overset{.}{\beta }\overset{..}{\theta }-\overset{.}%
{\theta }\overset{..}{\beta }\right) }{2\sqrt{f_{4}^{2}\overset{.}{\theta }^{2}%
-f_{2}^{2}\overset{.}{\beta }^{2}}}e_{3}+\frac{\left( f_{4}^{\prime }f_{2}%
\overset{.}{\beta }^{2}-f_{2}^{\prime }f_{4}\overset{.}{\theta }%
^{2}+f_{2}^{\prime \prime }f_{4}^{\prime }-f_{2}^{\prime }f_{4}^{\prime
\prime }\right) }{2\sqrt{f_{4}^{\prime 2}-f_{2}^{\prime 2}}}e_{4}
\end{equation*}
where
$e_{3}=\frac{\left( -f_{4}\overset{.}{\theta }\cos \beta ,f_{4}\overset%
{.}{\theta }\sin \beta ,-f_{2}\overset{.}{\beta }\cos \theta ,f_{2}\overset{.%
}{\beta }\sin \theta \right) }{\sqrt{-f_{2}^{2}\overset{.}{\beta }%
^{2}+f_{4}^{2}\overset{.}{\theta }^{2}}},e_{4}=\frac{\left( f_{4}^{\prime }\sin
\beta ,f_{4}^{\prime }\cos \beta ,f_{2}^{\prime }\sin \theta ,f_{2}^{\prime
}\cos \theta \right) }{\sqrt{-f_{2}^{\prime 2}+f_{4}^{\prime 2}}}$; $ -\infty
<x,y,z,\alpha ,\beta ,\theta <\infty ,s\in I$ and $f_{i}\in C^{\infty }$, \cite{5}.
\end{enumerate}
\end{theorem}

\section{Clairaut's theorem on the surfaces of rotation in $E_{2}^{4}$}

This section will use three different types of surfaces of rotation given
the previous section, and will generalize Clairaut's
theorem to these surfaces in $E_{2}^{4}$.

\subsection{Clairaut's theorem on the hyperbolic surface of rotation $%
\Upsilon^{1}$}

In this section, one will use the hyperbolic surface of rotation
parametrized as 
\begin{equation*}
\Upsilon^{1}(x,\alpha ,t)=\left( 
\begin{array}{c}
f_{1}\cosh x+f_{3}\sinh x,f_{2}\cosh \alpha +f_{4}\sinh \alpha , \\ 
f_{1}\sinh x+f_{3}\cosh x, f_{2}\sinh \alpha +f_{4}\cosh \alpha%
\end{array}%
\right).
\end{equation*}

Also, one can take the planar curve $\gamma $ for this surface of rotation
to be the intersection of $\Upsilon ^{1}(x,\alpha ,t)$ with $\varrho
,\vartheta =0($or $\xi ,\eta =0)$ for the coordinate system $(\xi ,\varrho
,\vartheta ,\eta ).$ Therefore, one can write that the curve $\gamma $ lies
on the $\xi \eta -$plane(or $\varrho \vartheta -$plane), and the curve can
be written by 
\begin{equation*}
\gamma (t)=(f_{1}(t),0,0,f_{4}(t));f_{1},f_{4}\in C^{\infty },
\end{equation*}%
then one has the parametrization%
\begin{equation*}
\Upsilon ^{1}(x,\alpha ,t)=\left( f_{1}\cosh x,f_{4}\sinh \alpha ,f_{1}\sinh
x,f_{4}\cosh \alpha \right) ;
\end{equation*}%
\begin{eqnarray*}
\Upsilon _{x}^{1}(x,\alpha ,t) &=&\left( f_{1}\sinh x,0,f_{1}\cosh
x,0\right) ;\Upsilon _{\alpha }^{1}(x,\alpha ,t)=\left( 0,f_{4}\cosh \alpha
,0,f_{4}\sinh \alpha \right) ; \\
\Upsilon _{t}^{1}(x,\alpha ,t) &=&\left( f_{1}^{\prime }\cosh
x,f_{4}^{\prime }\sinh \alpha ,f_{1}^{\prime }\sinh x,f_{4}^{\prime }\cosh
\alpha \right) .
\end{eqnarray*}

Hence, from the first fundamental form of the surface $\Upsilon^{1}$, one
has 
\begin{equation*}
\left\langle \Upsilon_{x}^{1},\Upsilon_{x}^{1}\right\rangle
=f_{1}^{2};\left\langle \Upsilon_{\alpha }^{1},\Upsilon_{\alpha
}^{1}\right\rangle =-f_{4}^{2};\left\langle
\Upsilon_{t}^{1},\Upsilon_{t}^{1}\right\rangle =f_{4}^{\prime
2}-f_{1}^{\prime 2};
\end{equation*}
\begin{equation}
\left\langle \Upsilon_{x}^{1},\Upsilon_{\alpha }^{1}\right\rangle
,\left\langle \Upsilon_{x}^{1},\Upsilon_{t}^{1}\right\rangle ,\left\langle
\Upsilon_{\alpha }^{1},\Upsilon_{t}^{1}\right\rangle =0;  \tag{3.1}
\end{equation}
\begin{equation*}
I_{\Upsilon^{1}}=%
\begin{pmatrix}
f_{1}^{2} & 0 & 0 \\ 
0 & -f_{4}^{2} & 0 \\ 
0 & 0 & f_{4}^{\prime 2}-f_{1}^{\prime 2}%
\end{pmatrix}%
\end{equation*}%
and one can write Lagrangian equation 
\begin{equation*}
L=f_{1}^{2}\overset{.}{x}^{2}-f_{4}^{2}\overset{.}{\alpha }^{2}+\left(
f_{4}^{\prime 2}-f_{1}^{\prime 2}\right) \overset{.}{t}^{2}.
\end{equation*}

So, the curve $\gamma $ is time-like and one writes $f_{4}^{\prime
2}-f_{1}^{\prime 2}=-1,f_{1}^{2}>0,-f_{4}^{2}$ $<0.$ Then, one gets 
\begin{equation*}
L=f_{1}^{2}\overset{.}{x}^{2}-f_{4}^{2}\overset{.}{\alpha }^{2}-\overset{.}{t%
}^{2}.
\end{equation*}

Hence, one can write following equations by using Clairaut's theorem,%
\begin{align*}
\frac{\partial }{\partial s}\left( \frac{\partial L}{\frac{\partial x}{%
\partial s}}\right) & =\frac{\partial L}{\partial x}\Rightarrow \frac{%
\partial }{\partial s}\left( 2f_{1}^{2}\overset{.}{x}\right) =0;\frac{%
\partial }{\partial s}\left( \frac{\partial L}{\frac{\partial \alpha }{%
\partial s}}\right) =\frac{\partial L}{\partial \alpha }\Rightarrow \frac{%
\partial }{\partial s}\left( -2f_{4}^{2}\overset{.}{\alpha }\right) =0, \\
\frac{\partial }{\partial s}\left( \frac{\partial L}{\frac{\partial t}{%
\partial s}}\right) & =\frac{\partial L}{\partial t}\Rightarrow \frac{%
\partial }{\partial s}\left( -2\overset{.}{t}\right) =2f_{1}^{\prime }f_{1}%
\overset{.}{x}^{2}-f_{4}^{\prime }f_{4}\overset{.}{\alpha }^{2}.
\end{align*}

Assume that $\gamma (t)$ is a geodesic on the surface $\Upsilon ^{1}.$
Hence, the curve $\gamma (t)$ can be written as 
\begin{equation*}
\overset{.}{\gamma }=\overset{.}{x}\Upsilon _{x}^{1}+\overset{.}{\alpha }%
\Upsilon _{\alpha }^{1}+\overset{.}{t}\Upsilon _{t}^{1},
\end{equation*}%
one can note that $\Upsilon _{x}^{1}=f_{1}N_{x}$ is a unit space-like vector
pointing along $x$-axis of the meridians, and $\Upsilon _{\alpha
}^{1}=f_{4}N_{\alpha }$ is a unit time-like vector pointing along the $%
\alpha $-axis of the parallels. Also $\Upsilon _{t}^{1}=N_{t}$ is a unit
time-like vector pointing along $t$-axis of the parallels. Also, the plane
spanned by $N_{t}$ and $N_{\alpha }$ are time-like and $\left\{ \Upsilon
_{x}^{1},\Upsilon _{\alpha }^{1},\Upsilon _{t}^{1}\right\} $ is an
orthonormal basis. Also, from (3.1), one gets 
\begin{equation}
\overset{.}{\gamma }=f_{1}\overset{.}{x}N_{x}+f_{4}\overset{.}{\alpha }%
N_{\alpha }+\overset{.}{t}N_{t}.  \tag{3.2}
\end{equation}

Note that if the $\gamma $ is time-like, since $N_{x}^{\bot }\in
Sp\{N_{\alpha },N_{t}\}$ one gets 
\begin{equation*}
\overset{.}{\gamma }=f_{1}N_{x}\overset{.}{x}+\left( f_{4}N_{\alpha }\overset%
{.}{\alpha }+\overset{.}{t}N_{t}\right) =N_{x}\cos \varphi _{1}+N_{x}^{\bot
}\sin \varphi _{1};
\end{equation*}%
\begin{equation}
=\cos \varphi _{1}N_{x}+\cosh \theta _{1}\sin \varphi _{1}N_{\alpha }+\sinh
\theta _{1}\sin \varphi _{1}N_{t},  \tag{3.3}
\end{equation}%
where $\varphi _{1}$ and $\theta _{1}$ are the angles between the meridians
of the surface and the time-like geodesic $\gamma .$

Also, from (3.2) and (3.3), one can write 
\begin{equation}
f_{1}\overset{.}{x}=\cos \varphi_{1} ;f_{4}\overset{.}{\alpha }=\cosh
\theta_{1} \sin \varphi_{1} ;\overset{.}{t}=\sinh \theta_{1} \sin
\varphi_{1} \text{.}  \tag{3.4}
\end{equation}

Similarly, for the curve $\gamma (t)=(0,f_{2}(t),f_{3}(t),0);f_{2},f_{3}\in
C^{\infty }.$ Then, the first fundamental form is written by%
\begin{equation*}
I_{S^{1}}=%
\begin{pmatrix}
-f_{3}^{2} & 0 & 0 \\ 
0 & f_{2}^{2} & 0 \\ 
0 & 0 & f_{3}^{\prime 2}-f_{2}^{\prime 2}%
\end{pmatrix}%
\end{equation*}%
and similar calculations are obtained.

Hence, from the Lagrangian equation, one has 
\begin{equation*}
2\overset{..}{t}=2f_{1}^{\prime }f_{1}\overset{.}{x}^{2}-f_{4}^{\prime }f_{4}%
\overset{.}{\alpha }^{2}.
\end{equation*}%
and 
\begin{equation*}
\frac{\partial }{\partial s}\left( 2f_{1}^{2}\overset{.}{x}\right)
=0\Rightarrow x=\frac{c_{1}t}{2f_{1}^{2}}+c_{3};\frac{\partial }{\partial s}%
\left( -2f_{4}^{2}\overset{.}{\alpha }\right) =0\Rightarrow \alpha =-\frac{%
c_{2}t}{2f_{4}^{2}}+c_{4},c_{i}\in 
\mathbb{R}
,
\end{equation*}%
which gives that $2f_{1}^{2}\overset{.}{x}$ and $-2f_{4}^{2}\overset{.}{%
\alpha }$ are constant along the geodesic curve. It follows that the
geodesics are given by 
\begin{equation*}
\frac{\partial }{\partial s}\left( 2f_{1}^{2}\overset{.}{x}\right) =0;\frac{%
\partial }{\partial s}\left( -2f_{4}^{2}\overset{.}{\alpha }\right) =0;\frac{%
\partial }{\partial s}\left( -2\overset{.}{t}\right) =2f_{1}^{\prime }f_{1}%
\overset{.}{x}^{2}-f_{4}^{\prime }f_{4}\overset{.}{\alpha }^{2}.
\end{equation*}

Now, by using the equations (3.4) and from the previous equations, one
writes 
\begin{align}
f_{1}\overset{.}{x}& =\cos \varphi_{1} \Rightarrow 2f_{1}^{2}\overset{.}{x}%
=2f_{1}\cos \varphi_{1} =cons.  \tag{3.5} \\
f_{4}\overset{.}{\alpha }& =\cosh \theta_{1} \sin \varphi_{1} \Rightarrow
-2f_{4}^{2}\overset{.}{\alpha }=-2f_{4}\cosh \theta_{1} \sin \varphi_{1}
=cons.  \tag{3.6} \\
-2\overset{.}{t}& =-2\sinh \theta_{1} \sin \varphi_{1} \neq cons.  \notag
\end{align}%
and for the equation $\frac{\partial }{\partial s}\left( \frac{\partial L}{%
\frac{\partial x}{\partial s}}\right) =\frac{\partial L}{\partial x},$ which
means that 
\begin{equation}
x=\int \frac{\cos \varphi_{1} }{f_{1}}ds  \tag{3.7}
\end{equation}%
is constant along the geodesic, conversely, if $\gamma $ is a curve with $%
2f_{1}\cos \varphi_{1}= $constant, the second equation is satisfied,
differentiating $L$ and substituting into the second Euler Lagrangian
equation yields the first Lagrangian equation. Furthermore, for the equation 
$\frac{\partial }{\partial s}\left( \frac{\partial L}{\frac{\partial \alpha 
}{\partial s}}\right) =\frac{\partial L}{\partial \alpha },$ 
\begin{equation}
\alpha =\int \frac{1}{f_{4}}\cosh \theta_{1} \sin \varphi_{1} ds  \tag{3.8}
\end{equation}%
is constant along the curve $\gamma .$ Therefore, one gives the following
theorem as a result of Clairaut's theorem on the
hyperbolic surface of rotation $\Upsilon^{1}\subset E_{2}^{4} $.

\begin{theorem}
Let $\gamma (t)=(f_{1}(t),0,0,f_{4}(t))$(or $\gamma
(t)=(0,f_{2}(t),f_{3}(t),0)$)$,f_{i}\in C^{\infty }$ be a time-like geodesic
curve on the hyperbolic surface of rotation $\Upsilon ^{1}$ in the $%
E_{2}^{4} $, let $f_{1}$ and $f_{4}$ be the distance functions from the axis
of rotation to a point on the surface. Therefore, $2f_{1}\cos \varphi _{1}$
and $-2f_{4}\cosh \theta _{1}\sin \varphi _{1}$ are constant along the curve 
$\gamma $ where $\varphi _{1}$ and $\theta _{1}$ are the angles between the
meridians of the surface and the time-like geodesic $\gamma $. Conversely,
if $2f_{1}\cos \varphi _{1}$ and $-2f_{4}\cosh \theta _{1}\sin \varphi _{1}$
are constant along $\gamma $, if no part of some parallels of the surface of
rotation, then $\gamma $ is time-like geodesic.
\end{theorem}

If one wants to obtain the general equation of geodesics, one should
consider the Euler-Lagrange equations%
\begin{align}
\overset{.}{x}=\frac{dx}{ds}=\frac{1}{f_{1}}\cos \varphi_{1};\overset{.}{%
\alpha }=\frac{d\alpha }{ds}=\frac{1}{f_{4}}\cosh \theta_{1} \sin\varphi_{1}
\tag{3.9}
\end{align}
then, adding (3.9) to Lagrangian equation $L,$ one has 
\begin{equation*}
L =f_{1}^{2}\left( \frac{dx}{ds}\right) ^{2}-f_{4}^{2}\left( \frac{d\alpha }{%
ds}\right) ^{2}-\left( \frac{dt}{dx}\frac{dx}{ds}\right) ^{2}
\end{equation*}
\begin{equation*}
\frac{dt}{dx}^{{}} = f_{1}\sqrt{1-\cosh ^{2}\theta_{1} \tan ^{2}\varphi_{1}
-L\sec ^{2}\varphi_{1} }
\end{equation*}%
or%
\begin{equation*}
L =f_{1}^{2}\left( \frac{dx}{ds}\right) ^{2}-f_{4}^{2}\left( \frac{d\alpha }{%
ds}\right) ^{2}-\left( \frac{dt}{d\alpha }\frac{d\alpha }{ds}\right) ^{2}
\end{equation*}
\begin{equation*}
\frac{dt}{d\alpha }^{{}}=f_{2}\sqrt{\cot ^{2}\varphi_{1} \tanh ^{2}
\theta_{1}-L\ sech ^{2} \varphi_{1} \ cosec ^{2} \varphi_{1}}.
\end{equation*}

\begin{theorem}
The general equation of geodesics on the hyperbolic surface of rotation $%
\Upsilon ^{1}(x,\alpha ,t)\subset E_{2}^{4}$, and for the parameters $%
\overset{.}{x}=\frac{1}{f_{1}}\cos \varphi _{1}$ and $\overset{.}{\alpha }=%
\frac{1}{f_{4}}\cosh \theta _{1}\sin \varphi _{1}$, are given by 
\begin{equation*}
\frac{dt}{dx}^{{}}=f_{1}\sqrt{1-\cosh ^{2}\theta _{1}\tan ^{2}\varphi
_{1}-L\sec ^{2}\varphi _{1}}
\end{equation*}%
or 
\begin{equation*}
\frac{dt}{d\alpha }^{{}}=f_{2}\sqrt{\cot ^{2}\varphi _{1}\tanh ^{2}\theta
_{1}-L\ sech^{2}\varphi _{1}\ cosec^{2}\varphi _{1}}.
\end{equation*}
\end{theorem}

\subsection{Clairaut's theorem on the hyperbolic surface of rotation $%
\Upsilon ^{2}(y,z,t)$}

In this section, one will use the hyperbolic surface of rotation
parametrized as 
\begin{equation*}
\Upsilon ^{2}(y,z,t)=\left( 
\begin{array}{c}
f_{1}\cosh y+f_{4}\sinh y,f_{2}\cosh z+f_{3}\sinh z, \\ 
f_{2}\sinh z+f_{3}\cosh z,f_{1}\sinh y+f_{4}\cosh y%
\end{array}%
\right) ,
\end{equation*}

Also, one can take the planar curve $\gamma $ for this surface of rotation
to be the intersection of $\Upsilon ^{2}(y,z,t)$ with $\vartheta ,\eta =0($%
or $\xi ,\varrho =0)$ for the coordinate system $(\xi ,\varrho ,\vartheta
,\eta )$, one can write that the curve $\gamma $ lies on the $\varrho \xi -$%
plane(or $\eta \vartheta -$plane), for $\gamma (t)=(f_{1}(t),f_{2}(t),0,0)$
or $\gamma (t)=(0,0,f_{3}(t),f_{4}(t));f_{i}\in C^{\infty },$ then one gets 
\begin{eqnarray*}
\Upsilon ^{2}(y,z,t) &=&\left( f_{1}\cosh y,f_{2}\cosh z,f_{2}\sinh
z,f_{1}\sinh y\right) \text{ or} \\
\Upsilon ^{2}(y,z,t) &=&\left( f_{4}\sinh y,f_{3}\sinh z,f_{3}\cosh
z,f_{4}\cosh y\right) .
\end{eqnarray*}

Here one will use the surface of rotation generated by the curve $\gamma
(t)=(f_{1}(t),f_{2}(t),0,0).$ Furthermore, 
\begin{eqnarray*}
\Upsilon_{y}^{2} &=&\left( f_{1}\sinh y,0,0,f_{1}\cosh y\right)
;\Upsilon_{z}^{2}=\left( 0,f_{2}\sinh z,f_{2}\cosh z,0\right) ; \\
\Upsilon_{t}^{2} &=&\left( f_{1}^{\prime }\cosh y,f_{2}^{\prime }\cosh
z,f_{2}^{\prime }\sinh x,f_{1}^{\prime }\sinh y\right),
\end{eqnarray*}
by resulting in the first fundamental form: 
\begin{equation*}
\left\langle \Upsilon_{y}^{2},\Upsilon_{y}^{2}\right\rangle
=f_{1}^{2};\left\langle \Upsilon_{z}^{2},\Upsilon_{z}^{2}\right\rangle
=f_{2}^{2};\left\langle \Upsilon_{t}^{2},\Upsilon_{t}^{2}\right\rangle
=-f_{2}^{\prime 2}-f_{1}^{\prime 2};
\end{equation*}
\begin{equation}
\left\langle \Upsilon_{y}^{2},\Upsilon_{z}^{2}\right\rangle ,\left\langle
\Upsilon_{y}^{2},\Upsilon_{t}^{2}\right\rangle ,\left\langle
\Upsilon_{z}^{2},\Upsilon_{t}^{2}\right\rangle =0;  \tag{3.10}
\end{equation}
\begin{equation*}
I_{\Upsilon^{2}}=%
\begin{pmatrix}
f_{1}^{2} & 0 & 0 \\ 
0 & f_{2}^{2} & 0 \\ 
0 & 0 & -f_{4}^{\prime 2}-f_{1}^{\prime 2}%
\end{pmatrix}%
\end{equation*}%
and one can write Lagrangian equation as follows 
\begin{equation*}
L=f_{1}^{2}\overset{.}{y}^{2}+f_{2}^{2}\overset{.}{z}^{2}+\left(
-f_{2}^{\prime 2}-f_{1}^{\prime 2}\right) \overset{.}{t}^{2}.
\end{equation*}

Here, one is interested the metric in $E_{2}^{4} $. So, one takes $\gamma $
to be time-like and one writes $-f_{2}^{\prime 2}-f_{1}^{\prime
2}=-1,f_{1}^{2},f_{2}^{2}>0$ $.$ Then, 
\begin{equation*}
L=f_{1}^{2}\overset{.}{y}^{2}+f_{2}^{2}\overset{.}{z}^{2}-\overset{.}{t}^{2}.
\end{equation*}

Hence, one can obtain the following equations using Clairaut's theorem, 
\begin{align*}
\frac{\partial }{\partial s}\left( \frac{\partial L}{\frac{\partial y}{%
\partial s}}\right) & =\frac{\partial L}{\partial y}\Rightarrow \frac{%
\partial }{\partial s}\left( 2f_{1}^{2}\overset{.}{y}\right) =0;\frac{%
\partial }{\partial s}\left( \frac{\partial L}{\frac{\partial z}{\partial s}}%
\right) =\frac{\partial L}{\partial z}\Rightarrow \frac{\partial }{\partial s%
}\left( 2f_{2}^{2}\overset{.}{z}\right) =0, \\
\frac{\partial }{\partial s}\left( \frac{\partial L}{\frac{\partial t}{%
\partial s}}\right) & =\frac{\partial L}{\partial t}\Rightarrow \frac{%
\partial }{\partial s}\left( 2\overset{.}{t}\right) =2f_{1}^{\prime }f_{1}%
\overset{.}{y}^{2}+2f_{2}^{\prime }f_{2}\overset{.}{z}^{2}.
\end{align*}

Let $\gamma (t)$ be a curve geodesic on the surface $\Upsilon^{2}.$ Hence, $%
\gamma (t)$ can be written as follows 
\begin{equation*}
\overset{.}{\gamma }=\overset{.}{y}\Upsilon_{y}^{2}+\overset{.}{z}%
\Upsilon_{z}^{2}+\overset{.}{t}\Upsilon_{t}^{2}.
\end{equation*}

So, one can note that $\Upsilon_{t}^{2}= N_{t} $ is a unit time-like vector
pointing along $t$-axis of the meridians, and $\Upsilon_{z}^{2}=f_{2}N_{z}$
is a unit space-like vector pointing along the $z$-axis of the parallels.
Also, $\Upsilon_{y}^{2}=f_{1}N_{y}$ is a unit space-like vector pointing
along $y$-axis of the parallels. It also follows that the plane spanned by $%
N_{y}$, $N_{z}$ is space-like and an orthonormal basis. Also, from (3.10),
one gets 
\begin{equation}
\overset{.}{\gamma }=\overset{.}{t}N_{t}+f_{1}\overset{.}{y}N_{y}+f_{2}%
\overset{.}{z}N_{z}.  \tag{3.11}
\end{equation}

Note that if $\gamma $ is time-like curve, since $N_{t}^{\bot }\in
Sp\{N_{y},N_{z}\}$, one gets 
\begin{equation*}
\overset{.}{\gamma }=f_{1}N_{y}\overset{.}{y}+f_{2}N_{z}\overset{.}{z}+%
\overset{.}{t}N_{t}=N_{t}\cosh \varphi _{2}+N_{t}^{\bot }\sinh \varphi _{2};
\end{equation*}%
\begin{equation}
=\cosh \varphi _{2}N_{t}+\cos \theta _{2}\sinh \varphi _{2}N_{y}+\sinh
\varphi _{2}\sin \theta _{2}N_{z},  \tag{3.12}
\end{equation}%
where $\varphi _{2}$ and $\theta _{2}$ are the angles between the meridians
of the surface and the time-like geodesic $\gamma ,$ and from (3.11) and
(3.12), one can write 
\begin{equation}
f_{1}\overset{.}{y}=\cos \theta _{2}\sinh \varphi _{2};f_{2}\overset{.}{z}%
=\sinh \varphi _{2}\sin \theta _{2};\overset{.}{t}=\cosh \varphi _{2}. 
\tag{3.13}
\end{equation}

Hence, from the Lagrangian equation one has 
\begin{equation*}
\overset{..}{t}=f_{1}^{\prime }f_{1}\overset{.}{y}^{2}+f_{2}^{\prime }f_{2}%
\overset{.}{z}^{2}
\end{equation*}%
and 
\begin{equation*}
\frac{\partial }{\partial s}\left( 2f_{1}^{2}\overset{.}{y}\right)
=0\Rightarrow y=\frac{c_{1}^{2}}{2f_{1}^{2}}t+c_{3}^{2};\frac{\partial }{%
\partial s}\left( 2f_{2}^{2}\overset{.}{z}\right) =0\Rightarrow z=\frac{%
c_{2}^{2}}{2f_{2}^{2}}t+c_{4}^{2},
\end{equation*}%
which gives that $2f_{1}^{2}\overset{.}{y}$ and $2f_{2}^{2}\overset{.}{z}$
are constant along the geodesic. It follows that the geodesics can be
written as 
\begin{equation*}
\frac{\partial }{\partial s}\left( 2f_{1}^{2}\overset{.}{y}\right) =0;\frac{%
\partial }{\partial s}\left( 2f_{2}^{2}\overset{.}{z}\right) =0;\frac{%
\partial }{\partial s}\left( -\overset{.}{t}\right) =f_{1}^{\prime }f_{1}%
\overset{.}{y}^{2}+f_{2}^{\prime }f_{2}\overset{.}{z}^{2}.
\end{equation*}

Now, from (3.13) one gets 
\begin{align}
f_{1}\overset{.}{y}& =\cos \theta_{2} \sinh \varphi_{2} \Rightarrow
2f_{1}^{2}\overset{.}{y}=2f_{1}\cos \theta_{2} \sinh \varphi_{2} =cons. 
\tag{3.14} \\
f_{2}\overset{.}{z}& =\sinh \varphi_{2} \sin \theta_{2} \Rightarrow
2f_{2}^{2}\overset{.}{z}=2f_{2}\sinh \varphi_{2} \sin \theta_{2} =cons. 
\tag{3.15} \\
2\overset{.}{t}& =2\cosh \varphi_{2} \neq cons.  \notag
\end{align}%
for the equation $\frac{\partial }{\partial s}\left( \frac{\partial L}{\frac{%
\partial y}{\partial s}}\right) =\frac{\partial L}{\partial y},$ 
\begin{equation}
y=\int \frac{\cos \theta_{2} \sinh \varphi_{2} }{f_{1}}ds  \tag{3.16}
\end{equation}%
is a constant, conversely, for the condition $2f_{1}\cos \theta_{2} \sinh
\varphi_{2}=$constant, the second equation is satisfied, differentiating $L$
and substituting into the second Euler Lagrangian equation yields the first
Lagrangian equation. Furthermore, for $\frac{\partial }{\partial s}\left( 
\frac{\partial L}{\frac{\partial z}{\partial s}}\right) =\frac{\partial L}{%
\partial z},$ 
\begin{equation}
z=\int \frac{\sinh \varphi_{2} \sin \theta_{2} }{f_{2}}ds  \tag{3.17}
\end{equation}%
is constant along the curve $\gamma .$ Hence, Clairaut's 
theorem is expressed on the hyperbolic surface of rotation given in $%
E_{2}^{4}$.

\begin{theorem}
Let $\gamma (t)=(f_{1}(t),f_{2}(t),0,0)$(or $\gamma
(t)=(0,0,f_{3}(t),f_{4}(t))$)$,f_{i}\in C^{\infty }$ be a time-like geodesic
curve on the hyperbolic surface of rotation $\Upsilon ^{2}$ in the $%
E_{2}^{4} $, and let $f_{1}$ and $f_{2}$ be the distance functions from the
axis of rotation to a point on the surface. Then, $2f_{1}\cos \theta
_{2}\sinh \varphi _{2}$ and $2f_{2}\sin \theta _{2}\sinh \varphi _{2}$ are
constant along the curve $\gamma $ where $\varphi _{2}$ and $\theta _{2}$
are the angles between the meridians of the surface and the time-like
geodesic curve $\gamma $. Conversely, if $2f_{1}\cos \theta _{2}\sinh
\varphi _{2}$ and $2f_{2}\sin \theta _{2}\sinh \varphi _{2}$ are constant
along the curve $\gamma $, if no part of some parallels of the surface of
rotation, then $\gamma $ is time-like geodesic.
\end{theorem}

In order to obtain the general equation of geodesics, one should consider
the Euler-Lagrange equations 
\begin{equation}
\overset{.}{y}=\frac{dy}{ds}=\frac{\cos \theta _{2}\sinh \varphi _{2}}{f_{1}}%
;\overset{.}{z}=\frac{dz}{ds}=\frac{\sinh \varphi _{2}\sin \theta _{2}}{f_{2}%
}.  \tag{3.18}
\end{equation}

By adding the equations (3.18) at Lagrangian equation $L,$ one has 
\begin{eqnarray*}
L &=&f_{1}^{2}\left( \frac{dy}{ds}\right) ^{2}+f_{2}^{2}\left( \frac{dz}{ds}
\right) ^{2}-\left( \frac{dt}{dy}\frac{dy}{ds}\right) ^{2} \\
\left( \frac{dt}{dx}\right) ^{2}\frac{\cos ^{2}\theta_{2} \sinh
^{2}\varphi_{2} }{f_{1}^{2}} &=&\sinh ^{2}\varphi_{2} -L\rightarrow \frac{dt%
}{dx}^{{}}=\frac{f_{1}\sqrt{\sinh ^{2}\varphi_{2} -L}}{\cos \theta_{2} \sinh
\varphi_{2} }
\end{eqnarray*}%
or%
\begin{eqnarray*}
L &=&f_{1}^{2}\left( \frac{dy}{ds}\right) ^{2}+f_{4}^{2}\left( \frac{dz}{ds}%
\right) ^{2}-\left( \frac{dt}{dz}\frac{dz}{ds}\right) ^{2} \\
\left( \frac{dt}{dz}\right) ^{2}\left( \frac{\sinh \varphi_{2} \sin
\theta_{2} }{f_{2}}\right) ^{2} &=&\sinh ^{2}\varphi_{2} -L\rightarrow \frac{%
dt}{dz}^{{}}=\frac{f_{2}\sqrt{\sinh ^{2}\varphi_{2} -L}}{\sinh \varphi_{2}
\sin \theta_{2} }.
\end{eqnarray*}

\begin{theorem}
The general equation of geodesics on the hyperbolic surface of rotation $%
\Upsilon ^{2}\subset E_{2}^{4}$, and for the parameters $\overset{.}{y}=%
\frac{\cos \theta _{2}\sinh \varphi _{2}}{f_{1}}$ and $\overset{.}{z}=\frac{%
\sinh \varphi _{2}\sin \theta _{2}}{f_{2}}$, are given by 
\begin{equation*}
\frac{dt}{dx}^{{}}=\frac{f_{1}\sqrt{\sinh ^{2}\varphi _{2}-L}}{\cos \theta
_{2}\sinh \varphi _{2}}\text{ or }\frac{dt}{dz}^{{}}=\frac{f_{2}\sqrt{\sinh
^{2}\varphi _{2}-L}}{\sinh \varphi _{2}\sin \theta _{2}}.
\end{equation*}
\end{theorem}

\subsection{Clairaut's theorem on the elliptic surfaces of rotation $%
\Upsilon ^{3}(\protect\beta ,\protect\theta ,s)$}

In this section, one will use the elliptic surface rotation parametrized as 
\begin{equation}
\Upsilon ^{3}(\beta ,\theta ,t)=\left( 
\begin{array}{c}
f_{1}\cos \beta +f_{2}\sin \beta ,-f_{1}\sin \beta +f_{2}\cos \beta , \\ 
f_{3}\cos \theta +f_{4}\sin \theta ,-f_{3}\sin \theta +f_{4}\cos \theta 
\end{array}%
\right) ,  \tag{3.19}
\end{equation}%
then one can take the planar curve $\gamma $ for this surface of rotation to
be the intersection of $\Upsilon ^{3}(\beta ,\theta ,s)$ with $\varrho ,\eta
=0($or $\xi ,\vartheta =0$) for the coordinate system $(\xi ,\varrho
,\vartheta ,\eta ).$ Therefore, the curve can be written by $\gamma
(t)=(0,f_{2}(t),0,f_{4}(t))$ (or $\gamma (t)=(f_{1}(t),0,f_{3}(t),0)$)$;$ $%
f_{i}\in C^{\infty },$ then one gets

\begin{equation*}
\Upsilon ^{3}(\beta ,\theta ,s)=\left( f_{2}\sin \beta ,f_{2}\cos \beta
,f_{4}\sin \theta ,f_{4}\cos \theta \right) \text{ }
\end{equation*}%
and resulting in the first fundamental form: 
\begin{eqnarray*}
\Upsilon _{\beta }^{3} &=&\left( f_{2}\cos \beta ,-f_{2}\sin \beta
,0,0\right) ;\Upsilon _{\upsilon }^{3}=\left( 0,0,f_{4}\cos \theta
,-f_{4}\sin \theta \right) ; \\
\Upsilon _{t}^{3} &=&\left( f_{2}^{\prime }\sin \beta ,f_{2}^{\prime }\cos
\beta ,f_{4}^{\prime }\sin \theta ,f_{4}^{\prime }\cos \theta \right) ;
\end{eqnarray*}%
\begin{equation*}
\left\langle \Upsilon _{\beta }^{3},\Upsilon _{\beta }^{3}\right\rangle
=-f_{2}^{2};\left\langle \Upsilon _{\theta }^{3},\Upsilon _{\theta
}^{3}\right\rangle =f_{4}^{2};\left\langle \Upsilon _{t}^{3},\Upsilon
_{t}^{3}\right\rangle =-f_{2}^{\prime 2}+f_{4}^{\prime 2};
\end{equation*}%
\begin{equation*}
\left\langle \Upsilon _{\beta }^{3},\Upsilon _{\theta }^{3}\right\rangle
,\left\langle \Upsilon _{\beta }^{3},\Upsilon _{t}^{3}\right\rangle
,\left\langle \Upsilon _{\theta }^{3},\Upsilon _{t}^{3}\right\rangle =0
\end{equation*}%
and 
\begin{equation*}
I_{\Upsilon ^{3}}=%
\begin{pmatrix}
-f_{2}^{2} & 0 & 0 \\ 
0 & f_{4}^{2} & 0 \\ 
0 & 0 & -f_{2}^{\prime 2}+f_{4}^{\prime 2}%
\end{pmatrix}%
.
\end{equation*}

Hence, one can write Lagrangian equation as follows 
\begin{equation*}
L=-f_{2}^{2}\overset{.}{\beta }^{2}+f_{4}^{2}\overset{.}{\theta }^{2}+\left(
-f_{2}^{\prime 2}+f_{4}^{\prime 2}\right) \overset{.}{t}^{2}.
\end{equation*}

If one takes $\gamma $ to be time-like, one can write $-f_{2}^{\prime
2}+f_{4}^{\prime 2}=-1,f_{4}^{2}>0$, $-f_{2}^{2}<0.$ Then, 
\begin{equation*}
L=-f_{2}^{2}\overset{.}{\beta }^{2}+f_{4}^{2}\overset{.}{\theta }^{2}-%
\overset{.}{t}^{2}.
\end{equation*}

Hence, one obtains the following equations using Clairaut\^{a}\euro
\texttrademark s theorem, 
\begin{align*}
\frac{\partial }{\partial s}\left( \frac{\partial L}{\frac{\partial \beta }{%
\partial s}}\right) & =\frac{\partial L}{\partial \beta }\Rightarrow \frac{%
\partial }{\partial s}\left( -2f_{2}^{2}\overset{.}{\beta }\right) =0;\frac{%
\partial }{\partial s}\left( \frac{\partial L}{\frac{\partial \theta }{%
\partial s}}\right) =\frac{\partial L}{\partial \theta }\Rightarrow \frac{%
\partial }{\partial s}\left( 2f_{4}^{2}\overset{.}{\theta }\right) =0, \\
\frac{\partial }{\partial s}\left( \frac{\partial L}{\frac{\partial t}{%
\partial s}}\right) & =\frac{\partial L}{\partial t}\Rightarrow \frac{%
\partial }{\partial s}\left( 2\overset{.}{t}\right) =-2f_{2}^{\prime }f_{2}%
\overset{.}{\beta }^{2}+2f_{4}^{\prime }f_{4}\overset{.}{\theta }^{2}.
\end{align*}

Let $\gamma $ be a geodesic curve on $\Upsilon ^{3}$. Hence, $\gamma (t)$
can be written as follows 
\begin{equation*}
\overset{.}{\gamma }=\overset{.}{\beta }\Upsilon _{\beta }^{2}+\overset{.}{%
\theta }\Upsilon _{\theta }^{2}+\overset{.}{t}\Upsilon _{t}^{2}.
\end{equation*}

Hence, one can note that $\Upsilon _{t}^{2}=N_{t}$ is a unit time-like
vector pointing along $t$-axis of the parallels, and $\Upsilon _{\beta
}^{2}=f_{2}N_{\beta }$ is a unit time-like vector pointing along the $\beta $%
-axis of the parallels. Also, $\Upsilon _{\theta }^{2}=f_{4}N_{\theta }$ is
a unit space-like vector pointing along $\theta $-axis of the meridians and
the plane spanned by $N_{t}$ and $N_{\beta }$ is time-like and an
orthonormal basis. Also, from (4.21), one gets 
\begin{equation}
\overset{.}{\gamma }=\overset{.}{t}N_{t}+f_{2}\overset{.}{\beta }N_{\beta
}+f_{4}\overset{.}{\theta }N_{\theta }.  \tag{3.20}
\end{equation}

Note that the $\gamma $ is time-like curve, since $N_{\theta }^{\bot }\in
Sp\{N_{t},N_{\beta }\}$, one gets 
\begin{equation*}
N_{\theta }^{\bot }=\cosh \theta _{3}N_{\beta }+\sinh \theta _{3}N_{t},
\end{equation*}%
\begin{equation}
\overset{.}{\gamma }=\overset{.}{t}N_{t}+f_{2}\overset{.}{\beta }N_{\beta
}+f_{4}\overset{.}{\theta }N_{\theta }=\cos \varphi _{3}N_{t}+\sin \varphi
_{3}\cosh \theta _{3}N_{\beta }+\sinh \theta _{3}\sin \varphi _{3}N_{\theta
}.  \tag{3.21}
\end{equation}%
where $\varphi _{3}$ and $\theta _{3}$ are the angles between the meridians
of the surface and the time-like geodesic curve $\gamma $.

Furthermore, from (3.20) and (3.21), one writes 
\begin{equation}
f_{2}\overset{.}{\beta }=\sin \varphi _{3}\cosh \theta _{3};f_{4}\overset{.}{%
\theta }=\sinh \theta _{3}\sin \varphi _{3};\overset{.}{t}=\cos \varphi _{3}.
\tag{3.22}
\end{equation}

Hence, from the Lagrangian equations, one has 
\begin{equation*}
-\overset{..}{t}=f_{2}^{\prime }f_{2}\overset{.}{\beta }^{2}+f_{4}^{\prime
}f_{4}\overset{.}{\theta }^{2};
\end{equation*}%
\begin{equation*}
\frac{\partial }{\partial s}\left( -2f_{2}^{2}\overset{.}{\beta }\right)
=0\Rightarrow \beta =\frac{-c_{1}^{3}t}{2f_{2}^{2}}+c_{3}^{3};\frac{\partial 
}{\partial s}\left( 2f_{4}^{2}\overset{.}{\upsilon }\right) =0\Rightarrow
\theta =\frac{c_{2}^{3}t}{2f_{4}^{2}}+c_{4}^{3},
\end{equation*}%
which gives that $2f_{2}^{2}\overset{.}{\beta }$ and $2f_{4}^{2}\overset{.}{%
\theta }$ are constant along the geodesic, and by using (3.22) one obtains 
\begin{align}
f_{2}\overset{.}{\beta }& =\sin \varphi _{3}\cosh \theta _{3}\Rightarrow
2f_{2}^{2}\overset{.}{\beta }=2f_{2}\sin \varphi _{3}\cosh \theta _{3}=cons.
\tag{3.23} \\
f_{4}\overset{.}{\theta }& =\sinh \theta _{3}\sin \varphi _{3}\Rightarrow
2f_{4}^{2}\overset{.}{\theta }=2f_{4}\sinh \theta _{3}\sin \varphi _{3}=cons.
\tag{3.24} \\
2\overset{.}{t}& =2\cos \varphi _{3}\neq cons.  \notag
\end{align}%
for $\frac{\partial }{\partial s}\left( \frac{\partial L}{\frac{\partial
\beta }{\partial s}}\right) =\frac{\partial L}{\partial \beta },$ one has $%
f_{2}\overset{.}{\beta }=\sin \varphi _{3}\cosh \theta _{3};f_{4}\overset{.}{%
\theta }=\sinh \theta _{3}\sin \varphi _{3};\overset{.}{t}=\cos \varphi _{3}$%
, 
\begin{equation}
\beta =\int \frac{\sin \varphi _{3}\cosh \theta _{3}}{f_{2}}ds  \tag{3.25}
\end{equation}%
is a constant, conversely for the condition $2f_{2}\sin \varphi _{3}\cosh
\theta _{3}=$constant and the equation $\frac{\partial }{\partial s}\left( 
\frac{\partial L}{\frac{\partial \theta }{\partial s}}\right) =\frac{%
\partial L}{\partial \theta },$ 
\begin{equation}
\upsilon =\int \frac{\sinh \theta _{3}\sin \varphi _{3}}{f_{4}}ds  \tag{3.26}
\end{equation}%
is constant along the curve $\gamma .$ Hence, the following theorem can be
given.

\begin{theorem}
Let $\gamma (t)=(0,f_{2}(t),0,f_{4}(t))$(or $\gamma
(t)=(f_{1}(t),0,f_{3}(t),0)$)$,f_{i}\in C^{\infty }$ be a time-like geodesic
curve on the elliptic surface of rotation $\Upsilon ^{3}\subset E_{2}^{4}$,
and let $f_{2}$ and $f_{4}$ be the distance functions from the axis of
rotation to a point on the surface. Then, $2f_{2}\sin \varphi _{3}\cosh
\theta _{3}$ and $2f_{4}\sinh \theta _{3}\sin \varphi _{3}$ are constant
along the curve $\gamma $ where $\varphi _{3}$ and $\theta _{3}$ are the
angles between the meridians of the surface and the time-like geodesic curve 
$\gamma $. Conversely, if $2f_{2}\sin \varphi _{3}\cosh \theta _{3}$ and $%
2f_{4}\sinh \theta _{3}\sin \varphi _{3}$ are constant along the curve $%
\gamma $, if no part of some parallels of the surface of rotation, then $%
\gamma $ is time-like geodesic curve.
\end{theorem}

For the general equation of geodesics, one should consider the
Euler-Lagrange equations 
\begin{equation}
\overset{.}{\beta }=\frac{d\beta }{ds}=\frac{\sin \varphi _{3}\cosh \theta
_{3}}{f_{2}};\overset{.}{\theta }=\frac{d\theta }{ds}=\frac{\sinh \theta
_{3}\sin \varphi _{3}}{f_{4}}.  \tag{3.27}
\end{equation}

By adding the equations in (3.27) at Lagrangian equation $L,$ one has 
\begin{equation*}
L=-f_{2}^{2}\left( \frac{d\beta }{ds}\right) ^{2}+f_{4}^{2}\left( \frac{%
d\theta }{ds}\right) ^{2}-\left( \frac{dt}{d\beta }\frac{d\beta }{ds}\right)
^{2}
\end{equation*}%
\begin{equation*}
\left( \frac{dt}{d\beta }\right) ^{2}\frac{\sin ^{2}\varphi _{3}\cosh
^{2}\theta _{3}}{f_{2}^{2}}=-\sin ^{2}\varphi _{3}-L;\frac{dt}{d\beta }%
^{{}}=i\frac{f_{2}\sqrt{L+\sin ^{2}\varphi _{3}}}{\sin \varphi _{3}\cosh
\theta _{3}}
\end{equation*}%
or%
\begin{equation*}
L=-f_{2}^{2}\left( \frac{d\beta }{ds}\right) ^{2}+f_{4}^{2}\left( \frac{%
d\theta }{ds}\right) ^{2}-\left( \frac{dt}{d\theta }\frac{d\theta }{ds}%
\right) ^{2}
\end{equation*}%
\begin{equation*}
\left( \frac{dt}{d\theta }\right) ^{2}\left( \frac{\sinh \theta _{3}\sin
\varphi _{3}}{f_{4}}\right) ^{2}=-\sin ^{2}\varphi _{3}-L\rightarrow \frac{dt%
}{d\theta }^{{}}=i\frac{f_{4}}{\sinh \theta _{3}\sin \varphi _{3}}\sqrt{\sin
^{2}\varphi _{3}+L}.
\end{equation*}

\begin{theorem}
The general equation of geodesics on the elliptic surface of rotation $%
\Upsilon ^{3}\subset E_{2}^{4}$, and for the parameters $\overset{.}{\beta }=%
\frac{\sin \varphi _{3}\cosh \theta _{3}}{f_{2}}$ and $\overset{.}{\theta }=%
\frac{\sinh \theta _{3}\sin \varphi _{3}}{f_{4}}$, are given by 
\begin{equation*}
\frac{dt}{d\beta }^{{}}=i\frac{f_{2}\sqrt{L+\sin ^{2}\varphi _{3}}}{\sin
\varphi _{3}\cosh \theta _{3}}\text{ \textit{or }}\frac{dt}{d\theta }^{{}}=i%
\frac{f_{4}}{\sinh \theta _{3}\sin \varphi _{3}}\sqrt{\sin ^{2}\varphi _{3}+L%
}.
\end{equation*}
\end{theorem}

\section{Conclusion}

This study generalizes Clairaut's theorem to pseudo Euclidean 4-space with
two index, and reviews Clairaut's theorem of surfaces of rotation which
define a well-known characterization of geodesics on a surface of rotation. 
Therefore, it is shown that the time-like geodesic curves on the hyperbolic
surface of rotation $\Upsilon ^{1}$ are completely characterized by $%
2f_{1}\cos \varphi _{1}$ and $-2f_{4}\cosh \theta _{1}\sin \varphi _{1}$
being constant, the time-like geodesics on the hyperbolic surface of
rotation $\Upsilon ^{2}$ are characterized by $2f_{1}\cos \theta _{2}\sinh
\varphi _{2}$ and $2f_{2}\sinh \varphi _{2}\sin \theta _{2}$ being constant,
and finally the time-like geodesics on the elliptic surface of rotation $%
\Upsilon ^{3}$ are characterized by $2f_{2}\sin \varphi _{3}\cosh \theta _{3}
$ and $2f_{4}\sinh \theta _{3}\sin \varphi _{3}$ being constant, where $%
\varphi _{i}$ and $\theta _{i}$ are the angles between the meridians of the
surface and the time-like geodesic curve $\gamma ,i=1,2,3$.

The authors are currently working on the properties of these surfaces of
rotation with a view to devising suitable metric in $E_{2}^{4}$ by adapting
the type of conservation laws considered in the paper. In our future
studies, the physical terms such as specific energy and specific angular
momentum will be examined with the help of the conditions obtained by using
the Clairaut's theorem for geodesics on these special surfaces.

\section{Acknowledgements}

The authors wish to express their thanks to the authors of literatures for
the supplied scientific aspects and idea for this study. Furthermore, we
request to explain great thanks to the reviewers for the constructive
comments and inputs given to improve the quality of our work.


\section{Funding}

Not applicable

\section{Conflict of Interest}

The authors declared that they have no conflict of interest.

\end{document}